\documentclass[twoside,notitlepage,12pt]{article}

\usepackage{amssymb}
\usepackage[leqno]{amsmath}
\usepackage{amsfonts}
\usepackage{latexsym}
\usepackage{amsopn}
\usepackage{amstext}
\usepackage{amsthm}

\usepackage[all]{xy}
\newdir{ >}{{}*!/-9pt/@{>}}

\usepackage{verbatim}
\usepackage[colorlinks]{hyperref}

\providecommand{\cal}{\mathcal}
\renewcommand{\Bbb}{\mathbb}

\newenvironment{pf}{\begin{proof}}{\end{proof}}

\newcommand{\Bee}{{\cal{B}}}

\newcommand{\Ef}{{\cal{F}}}

\newcommand{\El}{{\cal{L}}}

\newcommand{\Zee}{{\Bbb{Z}}}

\newcommand{\Nat}{{\Bbb{N}}}

\newcommand{\al}{\alpha}

\newcommand{\sig}{\sigma}

\renewcommand{\phi}{\varphi}
\renewcommand{\rho}{\varrho}

\newcommand{\rest}{\restriction}

\newcommand{\ntr}{{n\in\omega}}

\newcommand{\loe}{\leq}
\newcommand{\goe}{\geq}

\newcommand{\subs}{\subseteq}

\newcommand{\nnempty}{\ne\emptyset}

\newcommand{\diam}{\operatorname{diam}}

\newcommand{\Es}{{\cal{S}}}

\newcommand{\length}{\operatorname{length}}

\newtheorem{tw}{Theorem}[section]
\newtheorem{wn}[tw]{Corollary}

\newtheorem{claim}[tw]{Claim}

\newtheorem{twq}{Theorem}

\theoremstyle{definition}

\theoremstyle{remark}

\newcommand{\setof}[2]{\{#1\colon #2\}}

\newcommand{\sett}[2]{\{#1\}_{#2}}
\newcommand{\sn}[1]{\{#1\}} 
\newcommand{\map}[3]{#1\colon #2 \to #3} 
\newcommand{\img}[2]{#1[#2]} 

\newcommand{\dpower}[2]{[#1]^{#2}}

\providecommand{\nat}{\omega}

\newcommand{\Er}{\mathcal R}

\title{Perfect independent sets with respect to infinitely many relations}
\author{Martin Dole\v zal\thanks{Supported by RVO: 67985840.} \and Wies\l aw Kubi\'s\thanks{Supported by GA\v CR grant P201 14 07880S and RVO: 67985840.}\\
{\small Institute of Mathematics, Czech Academy of Sciences, Czech Republic}
}

\begin{document}

\maketitle

\begin{abstract}
We prove a result on perfect cliques with respect to countably many $G_{\delta}$ relations on a complete metric space.
As an application, we show that a Polish group contains a free subgroup generated by a perfect set as long as it contains any uncountable free subgroup.
This answers a recent question of G\l\c ab and Strobin.

\ \\


\end{abstract}

\tableofcontents

\section{Introduction}

Given an $n$-ary relation $R$ on a set $X$, a subset $S$ of $X$ is said to be \emph{$R$-independent} if for every sequence $s_1,\dots,s_n$ of pairwise distinct elements of $S$ it holds that
$$(s_1,\dots,s_n) \notin R.$$
This is one of the standard notions of independence considered widely in the literature, under various names.
The above concept becomes more useful after generalizing to arbitrary families of relations.
Namely, assuming $\Er$ is a family of relations on $X$ (possibly, each of different arity), we say that $S \subs X$ is \emph{$\Er$-independent} if it is $R$-independent for every $R \in \Er$.
Independence with respect to a family of relations was considered (with a slightly more technical definition) by Mycielski~\cite{MycielskiIndependent}.

Our aim is to present a dichotomy concerning the existence of independent sets in completely metrizable topological spaces.
Here is our main result:

\begin{tw}\label{ThmMejn}
	Let $X$ be a completely metrizable space of weight $\kappa \goe \aleph_0$, and let $\Er$ be a countable family of $F_\sigma$ relations on $X$.
	Then exactly one of the following two statements holds.
	\begin{enumerate}
		\item[\rm{(S)}] There exists an ordinal $\gamma < \kappa^+$ such that every $\Er$-independent set has Cantor-Bendixson rank $< \gamma$ (that is, $\gamma$th Cantor-Bendixson derivative of every $\Er$-in\-de\-pen\-dent set is empty).
		\item[\rm{(P)}] There exists a perfect $\Er$-independent set.
	\end{enumerate}
\end{tw}

Recall that a set $P$ is \emph{perfect} if it is nonempty, completely metrizable and has no isolated points.
Every perfect set contains a compact perfect set, namely, a topological copy of the Cantor set.
Note that the existence of a nonempty (even countable) dense-in-itself $\Er$-independent set already implies (P).
Furthermore, (P) holds whenever there exists an $\Er$-independent set of cardinality $>\kappa$, as such a set always contains a dense-in-itself subset of cardinality $>\kappa$.

Special case of Theorem~\ref{ThmMejn} has already been proved by the second author \cite{KubisGdelta}, under the assumptions that $X$ is separable and $\Er = \sn{R}$, where $R \subs X^n$ is a symmetric relation.

Theorem~\ref{ThmMejn} is parallel to the following result of Mycielski:

\begin{twq}[\cite{MycielskiIndependent}]
	Let $X$ be a completely metrizable space without isolated points, and let $\Er$ be a countable family of relations on $X$ such that each $R \in \Er$ is of the first category.
	Then there exists a perfect $\Er$-independent set.
\end{twq}
 
Clearly, every first category relation is contained in a first category  $F_\sig$ relation.
Thus, Mycielski's theorem says that if (S) holds in our dichotomy then at least one of the relations is not meager.

It is necessary to point out that the dichotomy totally fails when increasing the Borel complexity of the relations.
Namely, there exists a $G_\delta$ binary relation $R$ on the Cantor set $2^\omega$ such that every maximal $R$-independent set has cardinality exactly $\aleph_1$ and no perfect set is $R$-independent.
A concrete example of such a relation was found by Vejnar and the second author~\cite[Thm. 2.1]{KubisVejnar}, although its existence was proved earlier by Shelah~\cite{ShelahSquares}.

On the other hand, concerning a single binary relation, there is a significantly stronger result proved by Feng~\cite{Feng} for separable spaces and by Chaber and Pol~\cite{ChaberPol} for arbitrary spaces:

\begin{twq}\label{ThmquoteChaberPol}
	Let $X$ be a continuous image of a complete metric space of weight $\kappa\goe\aleph_0$, and let $R \subs X^2$ be a closed symmetric relation containing the diagonal of $X$.
	Then either $X = \bigcup_{\al<\kappa}X_\al$ such that each $X_\al^2 \subs R$ for every $\al < \kappa$ or there exists a perfect $R$-independent set.
\end{twq}

In fact, Chaber and Pol formulated and proved a more technical statement implying the one above.
Theorem~\ref{ThmquoteChaberPol} can be viewed as a far reaching generalization of the old and classical theorem of Suslin: every uncountable analytic set contains a perfect subset (just let $R$ be the diagonal of $X$).

A recent result of G\l\c ab and Strobin~\cite{GlaStr} asserts that a countable product of countable groups either contains a free group of cardinality $2^{\aleph_0}$ or else all of its free subgroups are countable.
They also ask whether such a dichotomy holds for all automorphism groups of countable first order structures.
In both cases, the groups carry a natural Polish topology, therefore our dichotomy provides an affirmative answer to the question of G\l\c ab and Strobin, providing additional information concerning the topological structure of sets of free generators.
The details are explained in Section~\ref{SecWolnoscAlg}.

Clearly, our main result (Theorem~\ref{ThmMejn}) can be rephrased dually, in terms of $G_\delta$ relations. 
Namely, given a relation $R \subs X^n$, a set $S \subs X$ is an \emph{$R$-clique}\footnote{The name \emph{clique} comes from graph theory, where $R$ is the edge relation.} if $(s_1,\dots,s_n) \in R$ whenever $s_1, \dots, s_n \in S$ are pairwise distinct.
Clearly, $S$ is an $R$-clique if and only if it is $R'$-independent, where $R' = X^n \setminus R$.
The notion of an \emph{$\Er$-clique}, where $\Er$ is a family of relations, is defined in the obvious way.
Now our dichotomy can be reformulated as follows:

\begin{tw}\label{ThmMejnGdelta}
	Let $X$ be a completely metrizable space of weight $\kappa \goe \aleph_0$, and let $\Er$ be a countable family of $G_\delta$ relations on $X$.
	Then exactly one of the following two statements holds.
	\begin{enumerate}
		\item[\rm{(S)}] There exists an ordinal $\gamma < \kappa^+$ such that every $\Er$-clique has Cantor-Bendixson rank $< \gamma$.
		\item[\rm{(P)}] There exists a perfect $\Er$-clique.
	\end{enumerate}
\end{tw}

The next section is devoted to the proof of this statement. Section~\ref{SecWolnoscAlg} explains applications to general topological algebras (including free groups), while the last section contains specific examples from group theory.

\section{Proof of the main result}\label{mainsection}

Recall that the \emph{weight} of a topological space is the least cardinality of an open base of the space.
A \emph{Polish space} is a completely metrizable separable topological space.
Recall that a topological space is \emph{analytic} if it is a continuous image of a Polish space.

For any subset $A$ of a topological space $X$ and for any ordinal $\gamma$, we denote by $A^{(\gamma)}$ the $\gamma$-th Cantor-Bendixson derivative of $A$. We also define the Cantor-Bendixson rank of $A$ as the least ordinal $\gamma$ such that $A^{(\gamma)}=\emptyset$. If such $\gamma$ does not exist, then the Cantor-Bendixson rank of $A$ is $+\infty$ which is, by definition, above all ordinals.


\begin{pf}[Proof of Theorem~\ref{ThmMejnGdelta}]
Let $\Er=\{R_n\}_{n\in\omega}$.
Then each $R_n$ is the intersection of countably many open relations $R_n^m$, $m\in\omega$.
If we put $\tilde\Er=\{R_n^m\}_{n,m\in\omega}$ then $\Er$-cliques and $\tilde\Er$-cliques coincide.
Therefore without loss of generality, we may suppose that $\Er$ consists of open relations.
	
Suppose that $X$ contains $\Er$-cliques of Cantor-Bendixson rank arbitrarily close to $\kappa^+$.
Let $\Er = \sett{R_n}{\ntr}$ and let $r_n$ be such that $R_n \subs X^{r_n}$.
Fix a complete metric $\rho \loe 1$ on $X$ and fix a base $\Bee$ of $X$ of cardinality $\loe\kappa$.
We construct inductively a Cantor tree $\sett{U_s}{s \in 2^{<\omega}} \subs \Bee$ with the following properties:
\begin{enumerate}
	\item[(1)] $\diam(U_s) \loe 2^{-\length(s)}$ for $s \in 2^{<\omega}$.
	\item[(2)] For every $k \in \nat$, for every $\gamma < \kappa^+$ there exists an $\Er$-clique $A_{k,\gamma}$ such that
$$A^{(\gamma)}_{k,\gamma} \cap U_s \nnempty$$
whenever $s \in 2^{k}$.
	\item[(3)]  For every $k \in \nat$, for every $i \loe k$, for every pairwise distinct $s_1, \dots, s_{r_i} \in 2^k$ it holds that
$$U_{s_1} \times \dots \times U_{s_{r_i}} \subs R_i.$$
\end{enumerate}
Note that (3) is automatically fulfilled when $r_i > 2^k$.
We start with $U_\emptyset := X$ if $r_0>1$ and with $U_\emptyset := R_0$ if $r_0=1$. Then (3) holds.
The existence of $\Er$-cliques of Cantor-Bendixson rank arbitrarily close to $\kappa^+$ guarantees (2).
Condition (1) holds by the choice of $\rho$.
Fix $n>0$ and suppose $\sett{U_s}{s \in 2^{<n}}$ has already been constructed.

\begin{claim}
Let $F$ be a finite $\Er$-clique and fix $\ntr$.
Then there exists a family $\sett{V_x}{x \in F}$ of pairwise disjoint open sets such that 
$x\in V_x$ for every $x\in F$ and such that for every $i \loe n$ and for every $S \in \dpower F {r_i}$ it holds that $$\prod_{x \in S}V_x \subs R_{i}.$$
\end{claim}

\begin{pf}
We may suppose that there is $i\loe n$ such that $r_i\loe|F|$, otherwise this is trivial. 
For each $i \loe n$ and for each $S \in \dpower F{r_i}$ choose a disjoint family of open sets $\sett{W^S_{x,i}}{x \in S}$ such that $x\in W^S_{x,i}$ for every $x\in S$ and such that
$$\prod_{x \in S}W^S_{x,i} \subs R_i.$$
Let
$$V_x = \bigcap \setof{W^S_{x,i}}{i \loe n, \; S \in \dpower F{r_i} \text{ is such that } x \in S}.$$
Then $\sett{V_x}{x\in F}$ is as required.
\end{pf}

Now, given $\gamma < \kappa^+$, choose $F_\gamma \subs A^{(\gamma)}_{n-1,\gamma+1} \cap \bigcup_{s \in 2^{n-1}}U_s$ such that $|F_\gamma \cap U_s| = 2$ for every $s \in 2^{n-1}$.
In particular, $|F_\gamma| = 2^n$.
Using the claim above, we find a disjoint family of open sets $\sett{V_{\gamma,x}}{x \in F_\gamma}$ such that
$x\in V_{\gamma,x}$ for every $x\in F_{\gamma}$ and such that
$$\prod_{x \in S}V_{\gamma,x} \subs R_i$$
whenever $i \loe n$ and $S \in \dpower F{r_i}$.
Without loss of generality, we may assume that the closure of each $V_{\gamma,x}$ is contained in the appropriate $U_s$ (where $s \in 2^{n-1}$ is such that $x \in U_s$) and that $V_{\gamma,x} \in \Bee$ for every $x \in F_\gamma$ and $\gamma < \kappa^+$.
We may also assume that each $V_{\gamma,x}$ has diameter $\loe 2^{-n}$.
Using the fact that $\Bee$ has cardinality $\loe\kappa$, there is a cofinal set $C \subs \kappa^+$ such that the family $\sett{V_{\gamma,x}}{x \in F_\gamma}$ is equal to some $\sett{V_s}{s \in 2^n}$ for every $\gamma \in C$, and $V_s \subs U_{s\rest (n-1)}$ for every $s \in 2^n$.
Define $U_s := V_s$ for $s \in 2^n$.
Then $\sett{U_s}{s \in 2^{\loe n}}$ clearly satisfies (1) and (3).
Condition (2) is witnessed by $A_{n, \gamma} := A_{n-1, \gamma+1}$ whenever $\gamma$ is from the cofinal set $C$.

Finally, the Cantor set resulting from the tree $\sett{U_s}{s \in 2^{<\omega}}$ is an $\Er$-clique, thanks to condition (3).
\end{pf}

We also note the following consequence of the main result.

\begin{tw}\label{Tweihrbiwgh2}
Let $X$ be a continuous image of a completely metrizable space of weight $\loe\kappa$ and let $\Er$ be a countable family of $G_{\delta}$ relations on $X$.
If there exists an $\Er$-clique of cardinality $\kappa^+$ then there also exists a perfect $\Er$-clique.
\end{tw}

\begin{pf}
	Fix a continuous surjection $\map f {\tilde X}X$, where $\tilde X$ is a completely metrizable space of weight $\loe\kappa$.
	Given $R \in \Er$ of arity $n$, define
	$$\tilde R = \setof{(x_1, \dots, x_n) \in \tilde X^n}{(f(x_1), \dots, f(x_n)) \in R}.$$
	Note that if $(x_1,\ldots,x_n)\in \tilde R$ then in particular $f\restriction_{\{x_1,\ldots,x_n\}}$ is one-to-one.
	The family $\tilde \Er = \sett{\tilde R}{R \in \Er}$ is clearly a countable family of $G_{\delta}$ relations on $\tilde X$.
	Suppose that $A\subseteq X$ is an $\Er$-clique of cardinality $\kappa^+$.
	Choose $\tilde A \subs f^{-1}(A)$ such that for every $a\in A$, it contains precisely one element from the preimage $f^{-1}(a)$.
	Then $\tilde A$ is an $\tilde \Er$-clique of cardinality $\kappa^+$.
	Knowing that $\tilde X$ has weight $\leq \kappa$, it is easy to see that for every $\gamma<\kappa^+$, we have $|\tilde A\setminus \tilde A^{(\gamma)}|\leq\kappa$, and so $\tilde A^{(\gamma)}\neq\emptyset$.
	It follows that the Cantor-Bendixson rank of $\tilde A$ is $\ge\kappa^+$.
	Finally, Theorem~\ref{ThmMejnGdelta} provides the existence of a perfect $\tilde \Er$-clique $P\subseteq \tilde X$.
	Assuming that $P$ is compact, $\img f P$ becomes a perfect $\Er$-clique in $X$.
\end{pf}

\begin{wn}\label{MainTheoremAnalytic}
Let $X$ be an analytic space and let $\Er$ be a countable family of $G_{\delta}$ relations on $X$.
If there exists an uncountable $\Er$-clique then there also exists a perfect $\Er$-clique.
\end{wn}

\section{Independence in topological algebras}\label{SecWolnoscAlg}

Recall that an \emph{abstract algebra} is a structure of the form $(X, \Ef)$, where $\Ef$ is a family of operations on the set $X$, where an \emph{operation} is simply a function $\map f {X^n}X$ with $n \goe1$ (the number $n$ is the \emph{arity} of $f$).
A \emph{topological algebra} is an algebra $(X,\Ef)$, where $X$ is a topological space and all operations in $\Ef$ are continuous.

In order to define independence in abstract algebras, we need to recall the notion of a term and equation.
Namely, assuming $\Ef = \sett{f^X_i}{i\in I}$, the \emph{language} of $(X,\Ef)$ is the set $\El = \sett{f_i}{i \in I}$ of formal operation symbols such that $f_i$ has the same arity as $f^X_i$.
A \emph{term} is, roughly speaking, an arbitrary formal operation that can be written as composition of formal operations from $\El$ and replacements of variables by other terms. The precise definition is recursive, of course. For details we refer to any textbook in model theory or logic.
Given two terms $t_1,t_2$ with variables $(x_1,\dots, x_m)$, the expression
$$t_1(x_1,\dots, x_m) = t_2(x_1,\dots,x_m)$$
is called an \emph{equation}.

Now let $\Es$ be a set of equations in $\El$ and let $(X,\Ef)$ be an algebra whose language is $\El$.
We say that $S \subs X$ is \emph{$\Es$-independent} if for every equation $t_1(x_1,\dots,x_m) = t_2(x_1,\dots,x_m)$ in $\Es$, it holds that
$$t_1(s_1,\dots, s_m) \ne t_2(s_1,\dots,s_m)$$
whenever $s_1, \dots, s_m \in S$ are pairwise distinct.
As an example, note that in the language of groups (where the operations are $(x_1,x_2) \mapsto x_1 \cdot x_2$, $x \mapsto x^{-1}$, and the constant $x \mapsto 1$) typical independence is with respect to all equations of the form $w(x_1,\dots, x_m) = 1$, where $w$ ranges over all nontrivial reduced words (which are of course terms).
Independent sets with respect to all these equations generate free subgroups.

Finally, notice that if $(X,\Ef)$ is a topological algebra whose underlying topology is Hausdorff, then for every equation $t_1(x_1,\dots,x_m) = t_2(x_1,\dots,x_m)$ the relation
$$R_{t_1,t_2} = \setof{(a_1,\dots,a_m) \in X^m}{t_1(a_1,\dots,a_m) = t_2(a_1,\dots,a_m)}$$
is closed.
Clearly, $R_{t_1,t_2}$-independence is the same as independence with respect to the equation $t_1 = t_2$.
Evidently, if the language $\El$ is countable then the set of all equations is countable, too.
Thus from our main result we obtain:

\begin{tw}\label{ThmFreeAlgsTops}
	Let $(X,\Ef)$ be a topological algebra of a countable language, whose underlying topology is metrizable and has weight $\kappa \goe \aleph_0$.
	Let $\Es$ be a fixed set of equations in the language of $(X,\Ef)$.
	Then exactly one of the following possibilities hold.
	\begin{enumerate}
		\item[{\rm (S)}] There exists an ordinal $\gamma < \kappa^+$ such that $\gamma$th Cantor-Bendixson derivative of every $\Es$-independent subset of $X$ is empty.
		\item[{\rm (P)}] There exists a perfect $\Es$-independent set.
	\end{enumerate}
	In particular, {\rm (P)} holds whenever there exists a dense-in-itself $\Es$-independent set.
\end{tw}

As a concrete example, we turn to the variety of groups. The main consequences of Theorem~\ref{ThmFreeAlgsTops} are the following.

\begin{wn}
Let $G$ be a completely metrizable topological group.
If $G$ has a dense-in-itself set of free generators then $G$ contains a free subgroup generated by a perfect set.
\end{wn}

\begin{wn}\label{CH}
Let $G$ be an analytic group. Then either all free subgroups of $G$ are countable or else $G$ contains a free subgroup generated by a perfect set.
\end{wn}

Note that the above corollary can be applied to any group of the form $\operatorname{Aut}(M)$ with $M$ a countable first order structure, as such a group carries a natural Polish topology (namely, it is a closed subgroup of $S_\infty$, the countable infinite permutation group).
This provides an affirmative answer to a question of G\l\c ab and Strobin~\cite{GlaStr}.

\section{Other examples from group theory}

In this section we demonstrate applications of our dichotomy to selected classes of groups.

\subsection{Free abelian subgroups}

It was proved in \cite[Theorem 4]{BalcerzykMycielski} that every locally compact non-$0$-dimensional group contains a free abelian subgroup generated by a set of cardinality $\mathfrak{c}$.
Here we provide examples of groups containing free abelian subgroups generated by a perfect set.

\begin{tw}
	\label{freeabelian}
	Let $G$ be a completely metrizable topological group of weight $\loe\kappa$.
	If for every ordinal $\gamma<\kappa^+$ there exists a set $S_\gamma \subs G$ of Cantor-Bendixson rank $\goe \gamma$ generating a free abelian group, then $G$ has a perfectly generated free abelian subgroup.
\end{tw}

\begin{pf}
	Let $R_n$ be the following $n$-ary relation: $(x_1, \dots, x_n) \in R_n$ if and only if
	$$x_1^{k_1}\ldots x_n^{k_n} \ne 1$$
	whenever $k_1,\ldots,k_n\in\Zee\setminus\{0\}$.
	Let $R$ be a binary relation defined by $(x_1,x_2)\in R$ if and only if $x_1,x_2,x_1^{-1},x_2^{-1}$ commute with each other.
	Clearly, each $R_n$ is a $G_{\delta}$ relation and $R$ is a closed relation.
	It is also clear that a subset of $G$ generates a free abelian group if and only if it is an $\Er$-clique, where $\Er = \{R_n\}_{n\in\Nat}\cup\{R\}$.
	Thus, the statement above follows from Theorem~\ref{ThmMejnGdelta}.
\end{pf}

\begin{wn}
	Let $G$ be a completely metrizable topological group.
	If $G$ has a free abelian subgroup which is dense-in-itself then $G$ contains a free abelian subgroup generated by a perfect set.
\end{wn}

\begin{wn}
	Let $G$ be an analytic group. Then either all free abelian subgroups of $G$ are countable or else $G$ contains a free abelian subgroup generated by a perfect set.
\end{wn}

\subsection{Torsion-free subgroups}

\begin{tw}\label{torsion-free}
	Let $G$ be a completely metrizable topological group of weight $\loe\kappa$.
	If for every ordinal $\gamma<\kappa^+$ there exists a set $S_\gamma \subs G$ of Cantor-Bendixson rank $\goe \gamma$ generating a torsion-free (resp. abelian torsion-free) group, then $G$ has a perfectly generated torsion-free (resp. abelian torsion-free) subgroup.
\end{tw}

\begin{pf}
	For $n,k\in\Nat$, let $R_{n,k}$ be the following $n$-ary relation: $(x_1, \dots, x_n) \in R_{n,k}$ if and only if either
	$$w(x_1,\ldots,x_n)=1$$
	or
	$$\left(w(x_1,\ldots,x_n)\right)^k \ne 1$$
	whenever $w$ is any word with the domain $G^n$.
	We put $\Er=\{R_{n,k}\}_{n,k\in\Nat}$.
	Clearly, each $R_{n,k}$ is a $G_{\delta}$ relation.
	It is also clear that a subset of $G$ generates a torsion-free group if and only if it is an $\Er$-clique.
	Thus, the statement for general torsion-free subgroups follows from Theorem~\ref{ThmMejnGdelta}.
	
	In the case of abelian torsion-free subgroups, we add one more binary relation $R$ defined by $(x_1,x_2)\in R$ if and only if $x_1,x_2,x_1^{-1},x_2^{-1}$ commute with each other. This relation is clearly closed (and therefore $G_{\delta}$).
\end{pf}

\begin{wn}
	Let $G$ be a completely metrizable topological group.
	If $G$ has a torsion-free (resp. abelian torsion-free) subgroup which is dense-in-itself then $G$ contains a torsion-free (resp. abelian torsion-free) subgroup generated by a perfect set.
\end{wn}

\begin{wn}
	Let $G$ be an analytic group. Then either all torsion-free (resp. abelian torsion-free) subgroups of $G$ are countable or else $G$ contains a torsion-free (resp. abelian torsion-free) subgroup generated by a perfect set.
\end{wn}

\subsection{Subgroups where all elements have a fixed order}

\begin{tw}\label{fixed_order}
	Let $F\subseteq \{2,3,\ldots\}$ be a finite set.
	Let $G$ be a completely metrizable topological group of weight $\loe\kappa$.
	If for every ordinal $\gamma<\kappa^+$ there exists a set $S_\gamma \subs G$ of Cantor-Bendixson rank $\goe \gamma$ generating a group where the order of each nonidentity element is in $F$, then $G$ has a perfectly generated subgroup where the order of each nonidentity element is in $F$.
\end{tw}

\begin{pf}
	For $n\in\Nat$, let $R_n$ be the following $n$-ary relation: $(x_1, \dots, x_n) \in R_n$ if and only for any word $w$ with the domain $G^n$ either $w(x_1,\ldots,x_n)=1$ or there is $k\in F$ such that
	\begin{itemize}
		\item[(i)] $\left(w(x_1,\ldots,x_n)\right)^i \ne 1$ for $i=1,\ldots,k-1$,
		\item[(ii)] $\left(w(x_1,\ldots,x_n)\right)^k = 1$.
	\end{itemize}
	We put $\Er=\{R_n\}_{n\in\Nat}$.
	Clearly, each $R_n$ is a $G_{\delta}$ relation.
	It is also clear that a subset of $G$ generates a group where the order of each nonidentity element is in $F$ if and only if it is an $\Er$-clique.
	Thus, the statement follows from Theorem~\ref{ThmMejnGdelta}.
\end{pf}

\begin{wn}
	Let $F\subseteq \{2,3,\ldots\}$ be a finite set.
	Let $G$ be a completely metrizable topological group.
	If $G$ has a dense-in-itself subgroup where the order of each nonidentity element is in $F$, then $G$ contains a perfectly generated subgroup where the order of each nonidentity element is in $F$.
\end{wn}

\begin{wn}
	Let $F\subseteq \{2,3,\ldots\}$ be a finite set and let $G$ be an analytic group. Then either all subgroups of $G$ where the order of each nonidentity element is in $F$ are countable, or else $G$ contains a perfectly generated subgroup where the order of each nonidentity element is in $F$.
\end{wn}

\subsection{CA-subgroups}

Recall that a group is said to be a \emph{CA-group} (or a \emph{centralizer abelian group}) if the centralizer of any nonidentity element is an abelian subgroup.
	
\begin{tw}\label{CA-group}
	Let $G$ be a completely metrizable topological group of weight $\loe\kappa$.
	If for every ordinal $\gamma<\kappa^+$ there exists a set $S_\gamma \subs G$ of Cantor-Bendixson rank $\goe \gamma$ generating a CA-subgroup then $G$ has a perfectly generated CA-subgroup.
\end{tw}
	
\begin{pf}
	For $n\in\Nat$, let $R_n$ be the following $n$-ary relation: $(x_1, \dots, x_n) \in R_n$ if and only for any words $w_i, i=1,2,3$, where the domain of each $w_i$ is $G^n$ it holds under the notation
	\begin{equation*}
		\tilde w_i=w_i(x_1,\ldots,x_n),\; i=1,2,3,
    \end{equation*}
	that
	\begin{equation*}
	    \tilde w_2=1\;\text{ or }\;
	    \tilde w_1\tilde w_2\ne\tilde w_2\tilde w_1
	    \;\text{ or }\;
	    \tilde w_3\tilde w_2\ne\tilde w_2\tilde w_3
	    \;\text{ or }\;
	    \tilde w_1\tilde w_3=\tilde w_3\tilde w_1.
	\end{equation*}
	We put $\Er=\{R_n\}_{n\in\Nat}$.
	Clearly, each $R_n$ is a $G_{\delta}$ relation.
	It is also clear that a subset of $G$ generates a CA-group if and only if it is an $\Er$-clique.
	Thus, the statement follows directly from Theorem~\ref{ThmMejnGdelta}.
\end{pf}
	
\begin{wn}
	Let $G$ be a completely metrizable topological group.
	If $G$ has a dense-in-itself CA-subgroup then $G$ contains a perfectly generated CA-subgroup.
\end{wn}
		
\begin{wn}
	Let $G$ be an analytic group. Then either all CA-subgroups of $G$ are countable, or else $G$ contains a perfectly generated CA-subgroup.
\end{wn}

\end{document}